\documentclass[12pt,a4paper]{amsart}

\usepackage[T1]{fontenc}
\usepackage[utf8]{inputenc}
\usepackage[british]{babel}
\usepackage{mathtools}
\usepackage{amsthm}
\usepackage{libertine}
\usepackage[libertine]{newtxmath}
\usepackage[mathscr]{euscript}
\usepackage{enumitem}
\usepackage{tikz-cd}
\usetikzlibrary{decorations.markings}
\usepackage{float}
\usepackage[
  backend=biber,
  style=alphabetic,
  maxnames=10,
  maxalphanames=10]{biblatex}
\usepackage{hyperref}
\usepackage[noabbrev]{cleveref}

\theoremstyle{plain}
\newtheorem{thm}{Theorem}
\newtheorem*{thm*}{Theorem}
\newtheorem{lm}[thm]{Lemma}
\newtheorem{prop}[thm]{Proposition}
\newtheorem{cor}[thm]{Corollary}

\theoremstyle{definition}
\newtheorem{defn}[thm]{Definition}
\newtheorem{construction}[thm]{Construction}
\newtheorem{nota}[thm]{Notation}
\newtheorem{exmp}[thm]{Example}

\theoremstyle{remark}
\newtheorem{rem}[thm]{Remark}
\Crefname{thm}{Theorem}{Theorems}
\Crefname{lm}{Lemma}{Lemmas}
\Crefname{prop}{Proposition}{Propositions}
\Crefname{cor}{Corollary}{Corollaries}
\Crefname{fact}{Fact}{Facts}
\Crefname{q}{Question}{Questions}
\Crefname{defn}{Definition}{Definitions}
\Crefname{construction}{Construction}{Constructions}
\Crefname{nota}{Notation}{Notations}
\Crefname{exmp}{Example}{Examples}
\Crefname{xca}{Exercise}{Exercises}
\Crefname{rem}{Remark}{Remarks}

\title[Adapted differentials as a qfh-sheaf]{Adapted differentials as a qfh-sheaf}
\author[Pedro N{\'u}{\~n}ez]{Pedro N{\'u}{\~n}ez}\thanks{The author gratefully acknowledges support by the DFG-Graduiertenkolleg GK1821 ``Cohomological Methods in Geometry'' at the University of Freiburg.}
\address{Pedro N{\'u}{\~n}ez \newline
\indent Albert-Ludwigs-Universit{\"a}t Freiburg, Mathematisches Institut \newline
\indent Ernst-Zermelo-Straße 1, 79104 Freiburg im Breisgau (Germany)}
\email{\normalfont\href{mailto:pedro.nunez@math.uni-freiburg.de}{pedro.nunez@math.uni-freiburg.de}}

\urladdr{\normalfont\href{https://home.mathematik.uni-freiburg.de/nunez/}{https://home.mathematik.uni-freiburg.de/nunez}}

\date{\today}

\makeatletter
\hypersetup{
  pdfauthor={\authors},
  pdftitle={\@title},
  colorlinks,
  linkcolor=[rgb]{0.2,0.2,0.6},
  citecolor=[rgb]{0.2,0.2,0.6},
  urlcolor=[rgb]{0.2,0.2,0.6}}
\makeatother

\begin{document}

\maketitle

\begin{abstract}
  We consider differential forms associated to Campana's geometric orbifolds from a new perspective, namely, as a qfh-sheaf on the variety underlying the geometric orbifold.
  This approach avoids having to choose a covering of the underlying variety, which is one of the drawbacks of a common way to work with these differential forms.
\end{abstract}

\tableofcontents

\section{Introduction}

\subsection{Motivation}

Geometric orbifolds play a central role in Campana's programme for the birational classification of algebraic varieties \cite{cam04}.
They are pairs $(X,\Delta)$ consisting of a normal variety $X$ and a Weil $\mathbb{Q}$-divisor of the form $\Delta = \sum_{i} \frac{m_{i} -  1}{m_{i}}D_{i}$ with $m_{i} \in \mathbb{N}_{\geq 2} \cup \{ \infty \}$ and with the convention that $\frac{\infty - 1}{\infty} = 1$.
These objects interpolate between projective varieties and logarithmic pairs in the sense of Iitaka \cite[\S 11]{iit82}.

As in the case of smooth projective varieties or logarithmic pairs, one of the fundamental tools to study the geometry of such an object are differential forms on it.
Associated to a geometric orbifold we should have differential forms with logarithmic poles of fractional order.
These differential forms are only well-defined on a covering space of our geometric orbifold, and this leads to the notion of adapted differentials.
This has two disadvantages:

\begin{itemize}
  \item The definition of adapted differentials requires the choice of a suitably ramified cover.
    This choice is not unique, and the properties of the resulting sheaf, such as the dimension of its space of global sections, depend on the chosen cover.
  \item These differentials do not form an object on $X$ a priori, like the sheaf of Kähler differentials on a projective variety or the sheaf of logarithmic differentials on a logarithmic pair do.
\end{itemize}

But there is a natural object on $X$ that takes all covers into account at once, namely, the presheaf on $\operatorname{Sch}/X$ described in \Cref{construction:presheaf}.
We show in \Cref{thm:main} that this presheaf is a qfh-sheaf on $\operatorname{Sch}/X$ in the sense of Suslin and Voevodsky.
The qfh-topology is recalled in \Cref{subsection:qfh}.

\subsection{Main result}

To make the following discussion a bit more self-contained, let us just say here that an \textit{adapted morphism} is a quasi-finite morphism $\gamma \colon Y \to X$ of normal varieties of the same dimension such that $\gamma^{*}\Delta$ has integer coefficients, and that on the domain of such a morphism we can define the sheaf of \textit{adapted differential $p$-forms} $\Omega_{(X,\Delta,\gamma)}^{[p]}$, which are the kind of differential forms that we want to consider.
See \Cref{defn:adapted}.

\begin{thm}\label{thm:main}
  Let $X$ be a normal variety over $\mathbb{C}$ and let $\Delta$ be a Weil $\mathbb{Q}$-divisor as above.
  For every $p \in \mathbb{N}$, there exists a presheaf $\Omega_{(X,\Delta)}^{p}$ on $\operatorname{Sch}/X$, unique up to isomorphism, with the following universal property:
  \begin{enumerate}[label=(\alph*)]
    \item\label{up:1} For every adapted morphism $\gamma \colon Y \to X$ there exists a morphism
      \[ \Omega_{(X,\Delta)}^{p}(Y) \to \Omega_{(X,\Delta,\gamma)}^{[p]}(Y), \]
      and these morphisms are compatible with the pull-back of differential forms.
    \item\label{up:2} For every presheaf $\mathscr{H}$ on $\operatorname{Sch}/X$ satisfying \ref{up:1}, there exists a unique morphism $\mathscr{H} \to \Omega_{(X,\Delta)}^{p}$ compatible with the morphisms from \ref{up:1}.
  \end{enumerate}
  Moreover, the morphisms in \ref{up:1} are isomorphisms.
  Furthermore, $\Omega_{(X,\Delta)}^{p}$ is a sheaf with respect to the qfh-topology on $\operatorname{Sch}/X$.
\end{thm}

\begin{rem}
  In particular, we can recover $\Omega_{(X,\Delta,\gamma)}^{[p]}$ as the restriction of $\Omega_{(X,\Delta)}^{p}$ to the small Zariski site on $Y$, because every dense open subset of the domain of an adapted morphism induces another adapted morphism.
\end{rem}

\begin{rem}
  Adapted differentials are in particular reflexive (logarithmic) differentials on the domain of an adapted morphism.
  Usual reflexive differentials on normal varieties are known to have qfh-descent, in the sense that \Cref{lm:galoisrestriction} also holds for reflexive differentials \cite[Théorème 2.2.2]{lw09}.
  However, reflexive differentials do not form a presheaf on the whole $\operatorname{Nor}/X$, because it is not always possible to define a meaningful pull-back, cf.~\cite[\S 1.2]{keb13} or \cite[\S B.1]{ks21}.
  Note also that the analogous qfh-descent property is known to fail for usual Kähler differentials \cite[Exemple 2.2.5]{lw09}.
\end{rem}

\begin{rem}
  The universal property in \Cref{thm:main} can be conveniently rephrased as a Kan extension, see \Cref{subsection:kan} and \Cref{rem:universalproperty}.
\end{rem}

Being a presheaf means that we can pull back sections along morphisms, and being a \emph{presheaf with transfers} means, roughly speaking, that we can more generally pull back sections along finite correspondences.

\begin{cor}
  \label{cor:transfers}
  In the setting of \Cref{thm:main}, the presheaf $\Omega_{(X,\Delta)}^{p}$ admits a natural structure of presheaf with transfers.
  \begin{proof}
    Follows from \Cref{thm:main} and \cite[Proposition 10.5.9]{cd19}.
  \end{proof}
\end{cor}

We refer to \cite[\S 10]{cd19} for a more detailed discussion on presheaves with transfers.

\subsection{Idea of the proof}\label{subsection:idea}

Adapted differentials form a Zariski sheaf on the domain of each adapted morphism, see \Cref{defn:adapted}.
Taking their global sections gives us a presheaf $\mathscr{F}$ on the full subcategory $\operatorname{Adpt}(X,\Delta) \subseteq \operatorname{Sch}/X$ of adapted morphisms, see \Cref{lm:presheaf}.
This presheaf can be extended to a presheaf $\Omega_{(X,\Delta)}^{p}$ on $\operatorname{Sch}/X$ as a Kan extension, see \Cref{construction:presheaf}.
The universal property and the description of the sections over adapted morphisms follow directly from the construction.

We can also express $\Omega_{(X,\Delta)}^{p}$ as the Kan extension of $\mathscr{G}$ to $\operatorname{Sch}/X$, where $\mathscr{G}$ is the Kan extension of $\mathscr{F}$ to the category $\operatorname{Nor}/X$ of integral normal schemes over $X$, cf.~\Cref{subsection:kan}.
By \cite[Lemma 6.4]{sv96}, it suffices to show that $\mathscr{G}$ satisfies the following properties:
\begin{enumerate}[label=(\roman*)]
  \item The presheaf $\mathscr{G}$ restricts to a Zariski sheaf on any integral normal scheme over $X$.\label{condition:1}
  \item For any integral normal scheme $T$ over $X$ and every finite Galois extension $L/\mathbb{C}(T)$, the restriction along the normalization $T' \to T$ of $T$ in $L$ gives a bijection $\mathscr{G}(T) \cong \mathscr{G}(T')^{\operatorname{Gal}(L/\mathbb{C}(T))}$, where the action of the Galois group on $\mathscr{G}(T')$ is given by pulling sections back along the corresponding automorphisms of $T'$ over $X$.\label{condition:2}
\end{enumerate}

Condition \ref{condition:1} is a formal consequence of the construction of the presheaf, cf.~\Cref{subsection:zariski}.
The idea for the surjectivity in condition \ref{condition:2} is the following.
Suppose $\pi \colon T' \to T$ is such a normalization of an integral normal scheme over $X$ in a finite Galois field extension with Galois group $G$ and suppose we are given some $\sigma' \in \mathscr{G}(T')^{G}$.
To define $\sigma \in \mathscr{G}(T)$ we need to specify an adapted differential $\sigma_{t}$ for every $t \colon Y \to T$ such that the composition $Y \to X$ is an adapted morphism, cf.~\Cref{construction:presheaf}.
We consider the commutative diagram
\begin{center}
  \begin{tikzcd}
    Y' \arrow[swap]{d}{t'} \arrow{r}{f} & Y \arrow{dr}{\gamma} \arrow{d}{t} & \\
    T' \arrow{r}{\pi} & T \arrow{r} & X
  \end{tikzcd}
\end{center}
in which $Y'$ is the normalization of the fiber product $T' \times_{T} Y$.
Assume for simplicity that $Y'$ is irreducible, so that the composition $Y' \to X$ is an adapted morphism.
We are forced by \Cref{construction:presheaf} to take an adapted differential $\sigma_{t}$ with the property that
\[ f^{*}(\sigma_{t}) = \sigma'_{t'}. \]
The group $G$ acts on $Y'$, and $(Y,f)$ is a quotient for this action in the sense of \cite[Exposé V, \S 1]{sga1}.
Therefore, the existence of such an adapted differential is ensured by the property that
\[ \Omega_{(X,\Delta,\gamma)}^{[p]}(Y) \cong \Omega_{(X,\Delta,\gamma \circ f)}^{[p]}(Y')^{G}, \]
which is proven in \Cref{lm:galoisrestriction}.

\subsection{Outline of this paper}

In \Cref{section:preliminaries} we discuss the necessary preliminaries on category theory and on adapted differentials.
Most of these preliminaries on adapted differentials are well-known to experts, but we include proofs for convenience.
In \Cref{section:proof} we prove \Cref{thm:main}.
In \Cref{section:computations} we compute the cohomology of $\Omega_{(X,\Delta)}^{p}$, as well as its sections over some quasi-finite schemes over $X$.

\subsection{Acknowledgements}

I would like to thank Stefan Kebekus for asking the question that motivated this result and for many useful discussions and comments on an earlier version of this paper.
I would like to thank Annette Huber-Klawitter for many useful discussions.
I would like to thank Shane Kelly for pointing out to me the relevance of presheaves with transfers.
I would also like to thank my fellow PhD students in Freiburg, especially Luca Terenzi, Vivien Vogelmann and Giovanni Zaccanelli, for the parallel study sessions which I spent thinking about this question.
Finally, I would like to thank an anonymous referee for many useful comments and suggestions.

\section{Preliminaries}\label{section:preliminaries}

\subsection{Notation and conventions}\label{subsection:notation}

We follow the notation and terminology in \cite{har77} and \cite{km98}.
In particular, varieties are assumed to be irreducible, and if $D$ is a Weil $\mathbb{Q}$-divisor on a normal variety $Y$, then we denote by $\lceil D \rceil$ (resp.~$\lfloor D \rfloor$) the divisor obtained by rounding up (resp.~down) its coefficients.
We work over the complex numbers $\mathbb{C}$ and fix a pair $(X,\Delta)$ as above once and for all.
Some other notation and conventions:
\begin{itemize}
  \item All schemes are assumed to be of finite type over $\mathbb{C}$.
  \item We denote by $\operatorname{Sch}/X$ the category of schemes over $X$.
  \item We denote by $\operatorname{Nor}/X$ the full subcategory of $\operatorname{Sch}/X$ consisting of integral normal schemes over $X$.
  \item A big open subset of a scheme is an open subset whose complement has codimension at least $2$.
  \item The quotient of a scheme by a group action will always mean a categorical quotient in the sense of \cite[Exposé V, \S 1]{sga1}.
  \item If $f \colon Z \to Y$ is a morphism of normal varieties and $\mathscr{A}$ is a coherent sheaf on $Y$, then we denote by $f^{[*]}\mathscr{A} := (f^{*}\mathscr{A})^{\vee\vee}$ the reflexive hull of $f^{*}\mathscr{A}$.
\end{itemize}
% TODO talk about the Kan extension here

\subsection{Reflexive differentials}\label{subsection:reflexive}

We follow \cite[\S 2.E]{gkkp11} when it comes to differential forms on normal varieties, i.e., we work with reflexive differential forms.
Let $(Y,D)$ be a pair consisting of a normal variety $Y$ and a Weil divisor $D$, and let $i \colon V \to Y$ be the inclusion of the largest open subset such that $(V,D|_{V})$ is snc.
Then we define the sheaf of \textit{reflexive logarithmic differential $p$-forms} as $\Omega_{Y}^{[p]}(\log{D}) := i_{*}\Omega_{V}^{p}(\log{D|_{V}})$ for every $p \in \mathbb{N}$.
Since $\Omega_{V}^{p}(\log{D|_{V}})$ is locally free and $V \subseteq Y$ is a big open subset, $\Omega_{Y}^{[p]}(\log{D})$ is indeed a reflexive coherent sheaf, and we may equivalently write it as the reflexive hull $\Omega_{Y}^{p}(\log{D})^{\vee\vee}$, see \cite[Notation 2.17]{gkkp11}.

In \Cref{subsection:pullback} we will say that \textit{the pull-back of Kähler differentials induces a morphism of reflexive logarithmic differentials} to mean that there is a morphism between sheaves of reflexive logarithmic differentials which agrees with the pull-back of Kähler differentials wherever this makes sense, e.g., over the complement of the support of the boundary divisor inside the smooth locus.

\subsection{Kan extensions}
\label{subsection:kan}
Kan extensions provide a convenient language to discuss a construction that will be used repeatedly in the proof of \Cref{thm:main}.
The geometric idea behind this construction is the following.
Suppose we are given a presheaf defined only on a basis of open subsets of a topological space.
Then, there is a natural way to extend it to a presheaf on the whole topological space.
Namely, by setting the sections over any given open subset to be the limit of the sections over basic open subsets contained in that open subset, cf.~\cite[Chap.~0, (3.2.1)]{egaI}.

In the categorical setting, let $u \colon \mathbf{A} \to \mathbf{B}$ be a functor.
Then we have a restriction functor on presheaves $u^{P} \colon \operatorname{PSh}(\mathbf{B}) \to \operatorname{PSh}(\mathbf{A})$, and this functor has a right adjoint ${_P}u \colon \operatorname{PSh}(\mathbf{A}) \to \operatorname{PSh}(\mathbf{B})$ which corresponds to the previous geometric construction, cf.~\cite[\href{https://stacks.math.columbia.edu/tag/00XF}{00XF}]{stacks-project}.
If $F$ is a presheaf on $\mathbf{A}$, then ${_P}u(F) = \operatorname{Ran}_{u}(F)$ is the \emph{right Kan extension of $F$ along $u$}, cf.~\cite[\S X.3]{mac98}.
We will use the following two observations:

\begin{rem}
  \label{rem:restrictionofextension}
  Let $i \colon \mathbf{A} \to \mathbf{B}$ be a fully faithful functor and let $F$ be a presheaf on $\mathbf{A}$.
  Then, the restriction of the Kan extension $\operatorname{Ran}_{i}(F)$ to $\mathbf{A}$ is canonically isomorphic to $F$ itself.
  See \cite[Corollary X.3.3]{mac98} and the proof of \Cref{lm:adaptedsections}.
\end{rem}

\begin{rem}
  \label{rem:compositionofextension}
  Let $u \colon \mathbf{A} \to \mathbf{B}$ and $v \colon \mathbf{B} \to \mathbf{C}$ be functors and let $F$ be a presheaf on $\mathbf{A}$.
  Since $(v \circ u)^{P} = u^{P} \circ v^{P}$, we have ${_P}(v \circ u) \cong {_P}v \circ {_P}u$.
  Therefore, $\operatorname{Ran}_{v \circ u}(F) \cong \operatorname{Ran}_{v}(\operatorname{Ran}_{u}(F))$.
\end{rem}

\subsection{qfh-topology}\label{subsection:qfh}

A morphism of schemes $p \colon Y \to Z$ is called a \emph{topological epimorphism} if the Zariski topology on $Z$ is the quotient topology of the Zariski topology on $Y$, i.e., if $p$ is surjective and a subset $A \subseteq Z$ is open if and only if $p^{-1}(A)$ is open.
It is called a \emph{universal topological epimorphism} if for any $W \to Z$ the projection from the fiber product $p_{W} \colon Y \times_{Z} W \to W$ is a topological epimorphism.
We recall the definition of the qfh-topology from \cite[Definition 3.1.2]{voe96}:

\begin{defn}[qfh-topology]
  The qfh-topology on the category of schemes is the Grothendieck topology associated with the pretopology whose coverings are of the form $\{ p_{i} \colon U_{i} \to Y \}$, where $\{ p_{i} \}$ is a finite family of quasi-finite morphisms of finite type such that the morphism $\sqcup p_{i} \colon \sqcup U_{i} \to Y$ is a universal topological epimorphism.
\end{defn}

\begin{rem}
  Let $Y$ be a normal connected scheme and let $\{ f_{i} \colon Y_{i} \to Y \}_{i \in I}$ be a finite family of quasi-finite morphisms of finite type such that all $Y_{i}$ are irreducible.
  Denote by $J$ the set of those $i$ for which $Y_{i}$ dominates $Y$.
  Then the above family is a qfh-covering if and only if $Y = \cup_{i \in J} f_{i}(Y_{i})$, cf.~\cite[Lemma 10.1]{sv96}.
\end{rem}

To prove \Cref{thm:main}, we will use the following:

\begin{lm}[{cf.~\cite[Lemma 6.4]{sv96}}]
  \label{lm:kan}
  Let $\mathscr{H}$ be a presheaf on $\operatorname{Nor}/X$.
  Assume that $\mathscr{H}$ satisfies the Zariski condition \ref{condition:1} and the Galois condition \ref{condition:2} spelled out in \Cref{subsection:idea}.
  Then the right Kan extension of $\mathscr{H}$ to $\operatorname{Sch}/X$ is a qfh-sheaf.
  \begin{proof}
    This right Kan extension is given by the right adjoint to the restriction functor $\operatorname{PSh}(\operatorname{Sch}/X) \to \operatorname{PSh}(\operatorname{Nor}/X)$, cf.~\Cref{rem:universalproperty}.
    Therefore, it is given by the functor $e \colon \operatorname{PSh}(\operatorname{Nor}/X) \to \operatorname{PSh}(\operatorname{Sch}/X)$ from \cite[\S 6]{sv96}.
    The arguments in the proof of \cite[Lemma 6.4]{sv96} do not use the assumption that we are working over a field, but only that we are working with finite type schemes over a field, cf.~\cite[Lemma 10.3]{sv96}.
    This is still true in our case, because $X$ is itself of finite type over $\mathbb{C}$, hence the claim.
  \end{proof}
\end{lm}

\subsection{Adapted morphisms and differentials}

Quasi-finite morphisms between normal varieties of the same dimension are open and dominant.
We can pull back Cartier divisors along dominant morphisms \cite[\href{https://stacks.math.columbia.edu/tag/02OO}{02OO}]{stacks-project}, and Weil divisors on a normal variety are Cartier over the regular locus, which is a big open subset.
The preimage of a big open subset under a quasi-finite morphism of varieties of the same dimension is again a big open subset.
Therefore, quasi-finite morphisms between normal varieties of the same dimension induce a pull-back homomorphism on Weil divisors.
This allows us to define:

\begin{defn}[Adapted morphism]\label{defn:adaptedmorphism}
  An \textit{adapted morphism} with respect to $(X,\Delta)$ is a quasi-finite morphism $\gamma \colon Y \to X$ such that $Y$ is a normal variety with $\dim(Y) = \dim(X)$ and such that $\gamma^{*}\Delta$ is a Weil divisor with integer coefficients.

  We denote by $\operatorname{Adpt}(X,\Delta)$ the full subcategory of $\operatorname{Sch}/X$ whose objects are adapted morphisms with respect to $(X,\Delta)$.
\end{defn}

\begin{exmp}
  Let $(X,\Delta)$ be a geometric orbifold with $\lceil \Delta \rceil - \Delta \neq 0$.
  Then $\operatorname{id}_{X}$ is a quasi-finite morphism between normal varieties of the same dimension, but it is not an adapted morphism with respect to $(X,\Delta)$.
\end{exmp}

\begin{exmp}
  \label{exmp:affineline}
  Let $X = \mathbb{A}^{1} = \operatorname{Spec}(\mathbb{C}[x])$ and let $\Delta = \frac{2}{3}\{ x = 0 \}$.
  Let $m \in \mathbb{N}$.
  The morphism $\gamma \colon \mathbb{A}^{1} \to \mathbb{A}^{1}$ given by $x \mapsto x^{m}$ is an adapted morphism with respect to $(X,\Delta)$ if and only if $m \in 3\mathbb{Z}$.
\end{exmp}

We will only consider a single pair $(X,\Delta)$, so from now we will simply refer to adapted morphisms with respect to $(X,\Delta)$ as adapted morphisms.

\begin{rem}\label{rem:composition}
  If $\gamma_{1} \colon Y_{1} \to X$ is an adapted morphism and $f \colon Y_{2} \to Y_{1}$ is a quasi-finite morphism between normal varieties of the same dimension, then $\gamma_{1} \circ f \colon Y_{2} \to X$ is again an adapted morphism.
\end{rem}

\begin{rem}\label{rem:cancellation}
  If $\gamma_{1} \colon Y_{1} \to X$ and $\gamma_{2} \colon Y_{2} \to X$ are adapted morphisms and $f \colon Y_{2} \to Y_{1}$ is a morphism over $X$, then $f$ is a quasi-finite morphism between normal varieties of the same dimension \cite[Proposition 12.17.(3)]{gw10}, hence open and dominant as well.
\end{rem}

We refer to \cite[II.\S 3]{del70} for the necessary background on logarithmic differentials and to \cite[\S 5]{cp19} for a discussion motivating the following:

\begin{defn}[Adapted differentials]\label{defn:adapted}
  Let $I$ be a set of indices such that $\Delta = \sum_{i \in I} \frac{m_{i} - 1}{m_{i}}D_{i}$ with $D_{i}$ distinct prime Weil divisors.
  We denote by $I_{0} \subseteq I$ the subset of indices such that $m_{i} \in \mathbb{N}$.
  Let $\gamma \colon Y \to X$ be an adapted morphism and assume first that $(X,\Delta)$ and $(Y,\gamma^{*}\Delta)$ are both snc pairs.
  Then $\gamma$ is flat, so $\gamma^{*}\mathscr{O}_{D_{i}} = \mathscr{O}_{\gamma^{*}D_{i}}$ for all $i \in I$, where $\gamma^{*}D_{i}$ is regarded here as a subscheme of $Y$ with the appropriate non-reduced structure.
  Moreover, since $\gamma^{*}D_{i} \geq \gamma^{*}(\frac{1}{m_{i}}D_{i})$ for all $i \in I_{0}$, we have a quotient morphism
  \[ q \colon \oplus_{i \in I} \gamma^{*}\mathscr{O}_{D_{i}} \to \oplus_{i \in I_{0}} \mathscr{O}_{\gamma^{*}(\frac{1}{m_{i}}D_{i})} \]
  which is an epimorphism of $\mathscr{O}_{Y}$-modules.
  We define the $\mathscr{O}_{Y}$-module of \textit{adapted differential $1$-forms} as
  \[ \Omega_{(X,\Delta,\gamma)}^{1} := \ker\left(\gamma^{*}\Omega_{X}^{1}(\log\lceil \Delta \rceil) \xrightarrow{\gamma^{*}(\operatorname{res})} \oplus_{i \in I} \gamma^{*} \mathscr{O}_{D_{i}} \xrightarrow{q} \oplus_{i \in I_{0}} \mathscr{O}_{\gamma^{*}(\frac{1}{m_{i}}D_{i})}\right). \]
  For each $p \in \mathbb{N}$, the $\mathscr{O}_{Y}$-module of \textit{adapted differential $p$-forms} is defined as
  \[ \Omega_{(X,\Delta,\gamma)}^{p} := \bigwedge^{p} \Omega_{(X,\Delta,\gamma)}^{1}, \]
  which is then locally free because $\Omega_{(X,\Delta,\gamma)}^{1}$ was locally free.

  In general, we consider the largest open subset $U \subseteq X$ such that $(X,\Delta)$ is snc over $U$.
  This is a big open subset, i.e., the codimension of its complement is at least $2$.
  Since $\gamma$ is quasi-finite between varieties of the same dimension, $\gamma^{-1}(U)$ is a big open subset as well.
  We can then consider the largest open subset $V \subseteq \gamma^{-1}(U)$ such that $(Y, \gamma^{*}\Delta)$ is snc over $V$, which is a big open subset of $Y$.
  We define the $\mathscr{O}_{Y}$-module of adapted differentials as
  \[ \Omega_{(X,\Delta,\gamma)}^{[p]} := i_{*}\left(\Omega_{(U,\Delta|_{U},\gamma|_{V})}^{p}\right), \]
  where $i \colon V \to Y$ denotes the open immersion.
\end{defn}

\begin{rem}\label{rem:reflexive}
  The sheaves $\Omega_{(X,\Delta,\gamma)}^{[p]}$ are coherent and reflexive, because they are the push-forward of locally free sheaves defined over a big open subset, see \cite[\href{https://stacks.math.columbia.edu/tag/0AY6}{0AY6}]{stacks-project}.
\end{rem}

\begin{exmp}
  If $\Delta = 0$, then any quasi-finite morphism $\gamma \colon Y \to X$ between normal varieties of the same dimension is an adapted morphism, and in this case we have $\gamma^{[*]}\Omega_{X}^{[p]} \cong \Omega_{(X,\Delta,\gamma)}^{[p]}$ for all $p \in \mathbb{N}$, cf.~also \Cref{lm:pullback} below.
\end{exmp}

\begin{exmp}
  In the setting of \Cref{exmp:affineline}, if $m = 3k$ for some $k \in \mathbb{N}_{>0}$, then $\Omega_{(X,\Delta,\gamma)}^{1}$ is the subsheaf of $\Omega_{\mathbb{A}^{1}}^{1}$ generated by $y^{k-1}dy$.
\end{exmp}

\begin{rem}
  \label{rem:rootstack}
  Root stacks have been used in the context of symmetric differentials on geometric orbifolds, see \cite[Proposition 3.4]{rou12}.
  If the construction of the (iterated) root stack of the boundary divisors is possible, e.g., if $X$ is locally factorial, then we can also understand adapted differentials in terms of root stacks.
  Namely, adapted differentials are the pull-back of differential forms defined on the root stack.
  To showcase the argument, let us look at a simple example.
  Let $X = \operatorname{Spec}(\mathbb{C}[x])$, $\Delta = \frac{m-1}{m}\{ x = 0 \}$ for some $m \in \mathbb{N}_{\geq 2}$, $Y = \operatorname{Spec}(\mathbb{C}[y])$ and $\gamma \colon Y \to X$ given by $x \mapsto y^{a}$ with $a = km$ for $k \in \mathbb{N}_{>0}$.
  Let $U = \operatorname{Spec}(\mathbb{C}[z]) \cong \operatorname{Spec}(\mathbb{C}[x,t]/(t^{m}-x))$ and consider the cyclic quotient $q \colon U \to X$ given by $x \mapsto z^{m}$.
  Denoting by $\mu_{m} \subseteq \mathbb{C}^{\times}$ the group of $m$-th roots of unity, we can describe the $m$-th root stack of $\{ x = 0 \} \subseteq X$ as the quotient stack $\mathcal{X} = [U/\mu_{m}]$, where the action of $\mu_{m}$ on $U$ is given by $\zeta \cdot z = \zeta z$ \cite[Theorem 10.3.10]{ols16}.
  The morphism $q$ factors as $\pi \circ u$, where $\pi \colon \mathcal{X} \to X$ is the coarse moduli space and $u \colon U \to \mathcal{X}$ is an atlas of the DM stack $\mathcal{X}$.
  Since $\gamma^{*}\{ x = 0 \} = km\{ y = 0 \}$ is an $m$-th power, the morphism $\gamma$ lifts to a morphism $\tilde{\gamma} \colon Y \to \mathcal{X}$ over $X$.
  This can in turn be lifted to the morphism $f \colon Y \to U$ given by $z \mapsto y^{k}$, and we obtain the following commutative diagram:
  \begin{center}
    \begin{tikzcd}
      & U \arrow{d}{u} \\
      Y \arrow{ur}{f} \arrow{r}{\tilde{\gamma}} \arrow[swap]{dr}{\gamma} & \mathcal{X} \arrow{d}{\pi} \\
      & X.
    \end{tikzcd}
  \end{center}
  In this situation, the pull-back of Kähler differentials induces an isomorphism $\tilde{\gamma}^{*}\Omega_{\mathcal{X}}^{1} \cong \Omega_{(X,\Delta,\gamma)}^{1}$.
  Indeed, by definition we have $u^{*}\Omega_{\mathcal{X}}^{1} = \Omega_{U}^{1}$, and the pull-back of Kähler differentials $f^{*}\Omega_{U}^{1} \to \Omega_{Y}^{1}$ induces an isomorphism onto the subsheaf generated by $y^{k-1}dy$, which is precisely $\Omega_{(X,\Delta,\gamma)}^{1}$.

  In more general situations, similar computations can be used to argue analytic locally over the smooth loci of $(X,\Delta)$ and $(Y,\gamma^{*}\Delta)$, cf.~proof of \Cref{lm:pullback}.
\end{rem}

\subsection{Pull-back of adapted differentials}\label{subsection:pullback}
From now on we fix a $p \in \mathbb{N}$.

The following lemma was observed by Stefan Kebekus and Erwan Rousseau in an upcoming work, but was probably known to experts before:

\begin{lm}\label{lm:pullback}
  Let $\gamma_{1} \colon Y_{1} \to X$ and $\gamma_{2} \colon Y_{2} \to X$ be two adapted morphisms and let $f \colon Y_{2} \to Y_{1}$ be a morphism over $X$.
  The pull-back of Kähler differentials induces an isomorphism
  \[ f^{[*]}\Omega_{(X,\Delta,\gamma_{1})}^{[p]} \cong \Omega_{(X,\Delta,\gamma_{2})}^{[p]}. \]
  \begin{proof}
    The restriction of sections of a reflexive coherent sheaf on a normal variety to a big open subset is an isomorphism \cite[Proposition 1.6]{har80}, so we may restrict to suitable big open subsets and assume that $(X,\Delta)$, $(Y_{1},\gamma_{1}^{*}\Delta)$ and $(Y_{2},\gamma_{2}^{*}\Delta)$ are all snc.
    Let $\Delta_{i}$ be the reduced divisor underlying $\gamma_{i}^{*}\lfloor \Delta \rfloor$ for each $i \in \{ 1, 2 \}$.
    Then $\Omega_{(X,\Delta,\gamma_{1})}^{p}$ is a locally free subsheaf of $\Omega_{Y_{1}}^{p}(\log{\Delta_{1}})$ and $f$ is dominant, so the pull-back of logarithmic differential forms as in \cite[\S 11.c]{iit82} induces an injective morphism
    \[ f^{*}\Omega_{(X,\Delta,\gamma_{1})}^{p} \to \Omega_{Y_{2}}^{p}(\log{\Delta_{2}}). \]
    The claim is that the image is precisely $\Omega_{(X,\Delta,\gamma_{2})}^{p}$.
    It suffices to show this on the stalks and for $p = 1$, so this can be checked analytic locally using the local generators in \cite[(60)]{cp19}.
    Indeed, let $z \in Y_{2}$ be a point and choose analytic local coordinates $x_{1}, \ldots, x_{n}$ around $\gamma_{2}(z) \in X$ such that
    \[ \Delta = \sum_{i = 1}^{r} \frac{m_{i} - 1}{m_{i}}D_{i} + \sum_{i = r+1}^{l}D_{i} \]
    with $D_{i} = \{ x_{i} = 0 \}$ around $\gamma_{2}(z)$.
    The simple normal crossing assumption implies that $l \leq n$, and for concreteness of notation let us assume that $1 \leq r < l < n$.
    We choose analytic local coordinates $y_{1}, \ldots, y_{n}$ around $f(z) \in Y_{1}$ such that
    \[ \gamma_{1}\left(y_{1}, \ldots, y_{n}\right) = \left( y_{1}^{a_{1}}, \ldots, y_{n}^{a_{n}} \right) \]
    with $a_{i} \in \mathbb{N}_{>0}$ for all $i \in \{1, \ldots, n\}$.
    Since $\gamma_{1}$ is an adapted morphism with respect to $(X,\Delta)$, for each $i \in \{1, \ldots, r \}$, there exists a $k_{i} \in \mathbb{N}$ such that $a_{i} = m_{i}k_{i}$.
    A system of local generators of $\Omega_{(X,\Delta,\gamma_{1})}^{1}$ around $f(z)$ is given by
    \[ y_{1}^{k_{1}-1}dy_{1}, \ldots, y_{r}^{k_{r}-1}dy_{r}, \frac{1}{y_{r+1}}dy_{r+1}, \ldots, y_{n}^{a_{n}-1}dy_{n}. \]
    This follows from the local descriptions of the morphisms involved in \Cref{defn:adapted}; cf.~also \cite[Definition 5.3]{cp19}, in which the adapted morphism would be chosen such that $k_{i} = 1$ for all $i \in \{1, \ldots, r\}$.
    Similarly, we choose analytic local coordinates $z_{1}, \ldots, z_{n}$ around $z \in Y_{2}$ such that
    \[ f\left(z_{1}, \ldots, z_{n}\right) = \left(z_{1}^{c_{1}}, \ldots, z_{n}^{c_{n}}\right) \]
    with $c_{i} \in \mathbb{N}_{>0}$ for all $i \in \{1, \ldots, n\}$.
    The morphism $\gamma_{2}$ is then given by
    \[ \gamma_{2}\left(z_{1}, \ldots, z_{n}\right) = \left(z_{1}^{c_{1}a_{1}}, \ldots, z_{n}^{c_{n}a_{n}}\right). \]
    The corresponding local generators of $\Omega_{(X,\Delta,\gamma_{2})}^{1}$ around $z$ are
    \[ z_{1}^{c_{1}k_{1}-1}dz_{1}, \ldots, z_{r}^{c_{r}k_{r}-1}dz_{r}, \frac{1}{z_{r+1}}dz_{r+1}, \ldots, z_{n}^{c_{n}a_{n}-1}dz_{n}, \]
    which are, up to multiplication by non-zero scalars, the pull-backs along $f$ of the local generators of $\Omega_{(X,\Delta,\gamma_{1})}^{1}$ around $f(z)$.
  \end{proof}
\end{lm}

\begin{lm}\label{lm:presheaf}
  If $\gamma_{1} \colon Y_{1} \to X$ and $\gamma_{2} \colon Y_{2} \to X$ are two adapted morphisms and $f \colon Y_{2} \to Y_{1}$ is a morphism over $X$, then the pull-back of Kähler differentials induces a $\mathbb{C}$-linear morphism $f^{*} \colon \Omega_{(X,\Delta,\gamma_{1})}^{[p]}(Y_{1}) \to \Omega_{(X,\Delta,\gamma_{2})}^{[p]}(Y_{2})$.
  Moreover:
  \begin{itemize}
    \item If $\gamma \colon Y \to X$ is an adapted morphism, then $\operatorname{id}_{Y}^{*} = \operatorname{id}_{\Omega_{(X,\Delta,\gamma)}^{[p]}(Y)}$.
    \item If $\gamma_{1} \colon Y_{1} \to X$, $\gamma_{2} \colon Y_{2} \to X$ and $\gamma_{3} \colon Y_{3} \to X$ are three adapted morphisms, and $f \colon Y_{2} \to Y_{1}$ and $g \colon Y_{3} \to Y_{2}$ are two morphisms over $X$, then $(f \circ g)^{*} = g^{*} \circ f^{*}$.
  \end{itemize}
  \begin{proof}
    The existence of $f^{*}$ follows from \cite[Proposition 1.6]{har80} and \Cref{lm:pullback}.
    Since the sheaves of adapted differentials are torsion-free, it suffices to check the desired identities over a dense open subset.
    But over a dense open subset these morphisms agree with the pull-back of Kähler differentials, hence the claim.
  \end{proof}
\end{lm}

\begin{nota}\label{nota:presheaf}
  We denote by $\mathscr{F}$ the presheaf on $\operatorname{Adpt}(X,\Delta)$ with $\mathscr{F}(Y) := \Omega_{(X,\Delta,\gamma)}^{[p]}(Y)$ for all adapted morphisms $\gamma \colon Y \to X$ and with restriction maps given as in \Cref{lm:presheaf}.
\end{nota}

\begin{lm}\label{lm:injectiverestriction}
  Let $\gamma_{1} \colon Y_{1} \to X$ and $\gamma_{2} \colon Y_{2} \to X$ be two adapted morphisms and let $f \colon Y_{2} \to Y_{1}$ be a morphism over $X$.
  Then the restriction map $f^{*} \colon \mathscr{F}(Y_{1}) \to \mathscr{F}(Y_{2})$ is injective.
  \begin{proof}
    Again by torsion freeness it suffices to check this over a dense open subset, over which it follows from the corresponding statement for pull-back of Kähler differentials along the dominant morphism $f$.
  \end{proof}
\end{lm}

\begin{lm}\label{lm:zrestriction}
  Let $\gamma \colon Y \to X$ be an adapted morphism and let $V \subseteq Y$ be an open subset.
  Then the restriction map $\mathscr{F}(Y) \to \mathscr{F}(V)$ is the usual restriction map of the sheaf $\Omega_{(X,\Delta,\gamma)}^{[p]}$.
  \begin{proof}
    Follows from torsion freeness and the corresponding statement for Kähler differentials as before.
  \end{proof}
\end{lm}

\begin{lm}\label{lm:galoisrestriction}
  Let $\gamma_{1} \colon Y_{1} \to X$ and $\gamma_{2} \colon Y_{2} \to X$ be two adapted morphisms and let $f \colon Y_{2} \to Y_{1}$ be a morphism over $X$ such that there exists a finite group $G$ acting on $Y_{2}$ in such a way that $(Y_{1},f)$ is a quotient of $Y_{2}$ by $G$.
  Then the restriction map $f^{*}$ induces a bijection $\mathscr{F}(Y_{1}) \to \mathscr{F}(Y_{2})^{G}$.
  \begin{proof}
    The $G$-action on $Y_{2}$ induces a $G$-sheaf structure on logarithmic differentials \cite[Fact 10.5]{gkkp11}, hence on adapted differentials by \Cref{lm:pullback}.
    The sheaf $(f_{*}\Omega_{(X,\Delta,\gamma_{2})}^{[p]})^{G}$ is then reflexive \cite[Lemma A.4]{gkkp11}, so we may restrict our attention to big open subsets and assume that $(X,\Delta)$, $(Y_{1},\gamma_{1}^{*}\Delta)$ and $(Y_{2},\gamma_{2}^{*}\Delta)$ are all snc.
    By \Cref{lm:pullback}, pull-back of Kähler differentials induces an isomorphism
    \[ f^{*}\Omega_{(X,\Delta,\gamma_{1})}^{p} \cong \Omega_{(X,\Delta,\gamma_{2})}^{p}. \]
    If $\phi_{g}$ denotes the automorphism of $Y_{2}$ corresponding to $g \in G$, then $f \circ \phi_{g} = f$, so the previous isomorphism is also an isomorphism of $G$-sheaves with respect to the $G$-sheaf structure induced on the pull-back.
    Hence, the assumption that $(Y_{1},f)$ is a quotient of $Y_{2}$ by $G$ and the projection formula imply that
    \[ \Omega_{(X,\Delta,\gamma_{1})}^{p} \cong \Omega_{(X,\Delta,\gamma_{1})}^{p} \otimes (f_{*}\mathscr{O}_{Y_{2}})^{G} \cong (f_{*}f^{*}\Omega_{(X,\Delta,\gamma_{1})}^{p})^{G} \cong (f_{*}\Omega_{(X,\Delta,\gamma_{2})}^{p})^{G}, \]
    and taking global sections yields the desired bijection.
  \end{proof}
\end{lm}

\section{Proof of \Cref{thm:main}}\label{section:proof}

To motivate \Cref{construction:presheaf}, we start with the following observation:

\begin{rem}\label{rem:universalproperty}
  If we denote by $i^{P} \colon \operatorname{PSh}(\operatorname{Sch}/X) \to \operatorname{PSh}(\operatorname{Adpt}(X,\Delta))$ the restriction functor, then the universal property in \Cref{thm:main} can be rephrased as follows:
  \begin{enumerate}[label=(\alph*)]
    \item There exists a morphism $\varepsilon \colon i^{P}(\Omega_{(X,\Delta)}^{p}) \to \mathscr{F}$.
    \item For every presheaf $\mathscr{H}$ on $\operatorname{Sch}/X$ such that there exists a morphism $\varphi \colon i^{P}(\mathscr{H}) \to \mathscr{F}$, there exists a unique morphism $\psi \colon \mathscr{H} \to \Omega_{(X,\Delta)}^{p}$ such that $\varphi = \varepsilon \circ i^{P}(\psi)$.
  \end{enumerate}
  In other words, $\Omega_{(X,\Delta)}^{p}$ is the right Kan extension of $\mathscr{F}$ along the inclusion $i \colon \operatorname{Adpt}(X,\Delta)^{\mathrm{op}} \to (\operatorname{Sch}/X)^{\mathrm{op}}$.
\end{rem}

\begin{construction}\label{construction:presheaf}
  We define $\Omega_{(X,\Delta)}^{p}$ as the right Kan extension of $\mathscr{F}$ to $\operatorname{Sch}/X$.
  Then $\Omega_{(X,\Delta)}^{p}$ satisfies the universal property in \Cref{thm:main} by construction, see \Cref{rem:universalproperty}.
  Recall from \Cref{subsection:kan} that we can write $\Omega_{(X,\Delta)}^{p} = {_P}i(\mathscr{F})$, where ${_P}i \colon \operatorname{PSh}(\operatorname{Adpt}(X,\Delta)) \to \operatorname{PSh}(\operatorname{Sch}/X)$ is the right adjoint to the restriction functor, as described in \cite[\href{https://stacks.math.columbia.edu/tag/00XF}{00XF}]{stacks-project}.
  Therefore, we can explicitly describe the sections of this presheaf as
  \begin{equation}\label{eqn:adjoint}
    \Omega_{(X,\Delta)}^{p}(T) = \lim_{\substack{Y \to T \\ Y \to X \text{ adapted}}} \mathscr{F}(Y).
  \end{equation}
  So a section $\sigma = (\sigma_{t})_{t} \in \Omega_{(X,\Delta)}^{p}(T)$ consists of a family of adapted differentials indexed by morphisms $t \colon Y \to T$ such that the composition $Y \to X$ is an adapted morphism.
  This family is moreover compatible in the sense that
  \begin{equation}\label{eqn:limit}
    f^{*}(\sigma_{t}) = \sigma_{t \circ f}
  \end{equation}
  for every morphism $f \colon Y' \to Y$ such that the composition $Y' \to X$ is an adapted morphism.
  If $\pi \colon T' \to T$ is a morphism and we are given a section $\sigma \in \Omega_{(X,\Delta)}^{p}(T)$, then $\Omega_{(X,\Delta)}^{p}(\pi)(\sigma)$ is the family uniquely determined by the property that
  \begin{equation}\label{eqn:restriction}
    \Omega_{(X,\Delta)}^{p}(\pi)(\sigma)_{t'} = \sigma_{\pi \circ t'}
  \end{equation}
  for every morphism $t' \colon Y \to T'$ such that the composition $Y \to X$ is an adapted morphism.
\end{construction}

\begin{lm}\label{lm:adaptedsections}
  Let $\gamma \colon Y \to X$ be an adapted morphism.
  Then the canonical morphism $\Omega_{(X,\Delta)}^{p}(Y) \to \mathscr{F}(Y)$ is an isomorphism.
  \begin{proof}
    This is a formal consequence of the inclusion $i \colon \operatorname{Adpt}(X,\Delta) \to \operatorname{Sch}/X$ being fully faithful, cf.~\Cref{rem:restrictionofextension}.
    Let us make this explicit in the case at hand.
    Since the inclusion $i \colon \operatorname{Adpt}(X,\Delta) \to \operatorname{Sch}/X$ is fully faithful, the indexing category in \Cref{eqn:adjoint} has an initial object, namely the identity on $Y$.
    Hence $\Omega_{(X,\Delta)}^{p}(Y) \cong \mathscr{F}(Y)$.
    More explicitly, the isomorphism is given by sending a family $(\sigma_{t})_{t}$ to the component $\sigma_{\operatorname{id_{Y}}}$ at $\operatorname{id}_{Y}$, and the inverse is given by pulling back along the corresponding morphisms.
  \end{proof}
\end{lm}

It remains to show that $\Omega_{(X,\Delta)}^{p}$ is a qfh-sheaf on $\operatorname{Sch}/X$.
As outlined in the introduction, we will apply \Cref{lm:kan} to the right Kan extension of $\mathscr{F}$ to $\operatorname{Nor}/X$.

\begin{nota}
  We denote by $\mathscr{G}$ the right Kan extension of $\mathscr{F}$ to $\operatorname{Nor}/X$, where $\mathscr{F}$ is the presheaf from \Cref{nota:presheaf}.
\end{nota}

By \Cref{rem:compositionofextension} and \Cref{rem:restrictionofextension}, we may also regard $\mathscr{G}$ as the restriction of $\Omega_{(X,\Delta)}^{p}$ to $\operatorname{Nor}/X$, so its sections are described by the same equations in \Cref{construction:presheaf} that describe the sections of $\Omega_{(X,\Delta)}^{p}$.

\subsection{Zariski condition \ref{condition:1}}
\label{subsection:zariski}
We need to show that for any integral normal $X$-scheme $T$, the presheaf $\mathscr{G}$ defines a sheaf in the small Zariski site on $T$ \cite[Definition 6.1]{sv96}.
Arbitrary fiber products do not exist in $\operatorname{Adpt}(X,\Delta)$ nor in $\operatorname{Nor}/X$, but fiber products along open immersions do exist.
Moreover, if $T \to X$ is in $\operatorname{Adpt}(X,\Delta)$ (resp.~in $\operatorname{Nor}/X$) and $U \to T$ is an open immersion, then the composition $U \to X$ is again in $\operatorname{Adpt}(X,\Delta)$ (resp.~in $\operatorname{Nor}/X$).
Therefore, taking Zariski coverings as coverings, both $\operatorname{Adpt}(X,\Delta)$ and $\operatorname{Nor}/X$ are sites in the sense of \cite[\href{https://stacks.math.columbia.edu/tag/00VH}{00VH}]{stacks-project}, and it suffices to show that $\mathscr{G}$ is a sheaf on $\operatorname{Nor}/X$ with respect to the corresponding Zariski topology.
Since the sheaves of adapted differentials are Zariski sheaves and the restriction maps of $\mathscr{F}$ are compatible with the restriction maps of these sheaves by \Cref{lm:zrestriction}, $\mathscr{F}$ is a sheaf with respect to the Zariski topology on $\operatorname{Adpt}(X,\Delta)$.
The inclusion functor $\operatorname{Adpt}(X,\Delta) \to \operatorname{Nor}/X$ is cocontinuous \cite[\href{https://stacks.math.columbia.edu/tag/00XJ}{00XJ}]{stacks-project}, so the sheaf property is preserved by the right Kan extension \cite[\href{https://stacks.math.columbia.edu/tag/00XK}{00XK}]{stacks-project} and $\mathscr{G}$ is a Zariski sheaf on $\operatorname{Nor}/X$.

\subsection{Galois condition \ref{condition:2}}
\label{subsection:galois}
In this subsection we fix a normal integral $X$-scheme $T \to X$ and a finite Galois extension $L/\mathbb{C}(T)$ of its function field $\mathbb{C}(T)$.
We denote by $G$ the Galois group of this extension and by $\pi \colon T' \to T$ the normalization of $T$ in $L$ \cite[Definition 12.42]{gw10}.

We will use the following observation.
This is not a new result, but we include a proof for convenience:

\begin{lm}\label{lm:galois}
  Let $Y$ and $Y'$ be normal integral schemes and let $f \colon Y' \to Y$ be a morphism.
  Then $f$ is a quotient of $Y'$ by the action of a finite group $H$ if and only if $f$ induces a finite Galois extension of function fields with Galois group $H$ and $Y'$ is the normalization of $Y$ in this extension of its function field.
  \begin{proof}
    Suppose $(Y,f)$ is a quotient of $Y'$ by a finite group $H$ in the sense of \cite[Exposé V, \S 1]{sga1}.
    Then $f$ is finite and the statement that we want to show is local on $Y$, so we may assume that both $Y'$ and $Y$ are affine, say $Y' = \operatorname{Spec}(B)$ and $Y = \operatorname{Spec}(A)$.
    The group $H$ acts on $B$ and the quotient morphism $f$ allows us to identify $A = B^{H}$, cf.~\cite[Exposé V, Proposition 1.1]{sga1}.
    By \cite[Exercise 5.12]{am69} we have an induced $H$-action on $B_{(0)}$ such that $A_{(0)} = (B_{(0)})^{H}$, so $\mathbb{C}(Y) \subseteq \mathbb{C}(Y')$ is a finite Galois extension with Galois group $H$.
    It remains to show that $Y'$ is the normalization of $Y$ inside $\mathbb{C}(Y')$, i.e., we want to show that $B = \{ b \in B_{(0)} \mid b \text{ integral over } A \}$.
    Every element in $B$ is integral over $A$, because $A \subseteq B$ is integral.
    And conversely, if $b \in B_{(0)}$ is integral over $A$, then it is also integral over $B$.
    Since $Y'$ is normal, we deduce that $b \in B$.

    Conversely, suppose that $Y'$ is the normalization of $Y$ in $\mathbb{C}(Y')$ and that $\mathbb{C}(Y) \subseteq \mathbb{C}(Y')$ is a finite Galois extension with Galois group $H$.
    We show that $H$ acts on $Y'$ with quotient $(Y,f)$.
    The statement is local on $Y$ and normalization is a finite morphism in our setting, so we may assume again that both $Y' = \operatorname{Spec}(B)$ and $Y = \operatorname{Spec}(A)$ are affine.
    By assumption there is an $H$-action on $B_{(0)}$ such that $A_{(0)} = (B_{(0)})^{H}$, and we want to show that it induces an $H$-action on $B$ such that $A = B^{H}$.
    Given $b \in B$ non-zero and $h \in H$, since $A \subseteq B$ is an integral extension, there exists a monic polynomial $P \in A[T]$ such that $P(b) = 0$.
    Since the coefficients of $P$ are $G$-invariant, this implies that $P(h \cdot b) = 0$, so $h\cdot b$ is also integral over $A$, hence over $B$.
    Since $Y'$ is normal, $h\cdot b \in B$.
    So the $H$-action on $B_{(0)}$ induces an $H$-action on $B$.
    Since $A_{(0)} = (B_{(0)})^{H}$ and $Y$ is normal, we have $A = B^{H}$.
  \end{proof}
\end{lm}

So in our current setting $\pi \colon T' \to T$ is the quotient of $T'$ by the induced $G$-action on $T'$.
For each $g \in G$ we denote by $\phi_{g}$ the corresponding automorphism of $T'$.
Then we have a $G$-action on $\mathscr{G}(T')$ given by pulling back along the corresponding automorphism.
The restriction map $\mathscr{G}(\pi) \colon \mathscr{G}(T) \to \mathscr{G}(T')$ has image contained in $\mathscr{G}(T')^{G}$, because  $\pi \circ \phi_{g} = \pi$ for all $g \in G$.
Our goal in this subsection is to show that this restriction map is a bijection onto $\mathscr{G}(T')^{G}$.

\begin{lm}\label{lm:injective}
  The morphism $\mathscr{G}(\pi) \colon \mathscr{G}(T) \to \mathscr{G}(T')^{G}$ is injective.
  \begin{proof}
    Let $\sigma, \sigma' \in \mathscr{G}(T)$ be two sections such that $\mathscr{G}(\pi)(\sigma) = \mathscr{G}(\pi)(\sigma')$.
    We want to show that $\sigma = \sigma'$, so let $t \colon Y \to T$ be a morphism such that the composition $Y \to X$ is an adapted morphism and let us show that $\sigma_{t} = \sigma'_{t}$.
    Let $Y'$ be a component of the normalization of the fiber product $T' \times_{T} Y$ that surjects onto $Y$, considered as a subscheme with its induced reduced structure, and let $f \colon Y' \to Y$ be the induced surjective morphism, which is also finite because so are the normalization, the closed immersion and the pull-back of $\pi$.
    Let $t' \colon Y' \to T'$ be the induced morphism, so that $\pi \circ t' = t \circ f$ and the composition $Y' \to X$ is an adapted morphism as well.
    It follows then from \Cref{eqn:limit} and \Cref{eqn:restriction} that $f^{*}(\sigma_{t}) = \sigma_{\pi \circ t'} = \sigma'_{\pi \circ t'} = f^{*}(\sigma'_{t})$, hence $\sigma_{t} = \sigma'_{t}$ by \Cref{lm:injectiverestriction}.
  \end{proof}
\end{lm}

\begin{lm}\label{lm:surjective}
  The morphism $\mathscr{G}(\pi) \colon \mathscr{G}(T) \to \mathscr{G}(T')^{G}$ is surjective.
  \begin{proof}
    Let $\sigma' \in \mathscr{G}(T')$ be a $G$-invariant section, i.e., $\mathscr{G}(\phi_{g})(\sigma') = \sigma'$ for all $g \in G$, where $\phi_{g}$ denotes the automorphism of $T'$ corresponding to $g$.
    We want to find a section $\sigma \in \mathscr{G}(T)$ such that $\mathscr{G}(\pi)(\sigma) = \sigma'$, so let $t \colon Y \to T$ be a morphism such that the composition $Y \to X$ is an adapted morphism.
    Let $Y_{0}'$ be any irreducible component of the fiber product $T' \times_{T} Y$, considered as a subscheme with its induced reduced structure.
    Let $G_{0}$ be the stabilizer of $Y_{0}'$ in $G$, so that the induced morphism $Y_{0}' \to Y$ is a quotient of $Y_{0}'$ by $G_{0}$, cf.~\cite[Lemma 2.1.11]{hub00} and \cite[Corollary 5.10]{sv96}.
    Let $\bar{Y}_{0}'$ be the normalization of $Y_{0}'$.
    Then $G_{0}$ acts on $\bar{Y}_{0}'$ as well and the induced morphism $f^{Y}_{0} \colon \bar{Y}_{0}' \to Y$ is still a quotient for this action, cf.~\Cref{lm:galois}.
    Let $t_{0} \colon \bar{Y}_{0}' \to T'$ be the induced morphism, so that $t \circ f^{Y}_{0} = \pi \circ t_{0}$ and the composition $\bar{Y}_{0}' \to X$ is an adapted morphism:
    \begin{center}
      \begin{tikzcd}
        \bar{Y}_{0}' \arrow[swap, bend right=45]{ddd}{t_{0}} \arrow{d}{\text{norm.}} \arrow{r}{f^{Y}_{0}} & Y \arrow[equal]{dd} \\
        Y'_{0} \arrow[hook]{d} & \\
        T' \times_{T} Y \arrow{d} \arrow{r} & Y \arrow{d}{t} \\
        T' \arrow{r}{\pi} & T.
      \end{tikzcd}
    \end{center}
    The desired $\sigma \in \mathscr{G}(T)$ must then have the property that $(f^{Y}_{0})^{*}(\sigma_{t}) = \sigma_{\pi \circ t_{0}} = \sigma'_{t_{0}}$.
    Let now $g \in G_{0}$ and let $\psi_{g}$ denote the corresponding automorphism of $\bar{Y}_{0}'$.
    From \Cref{eqn:limit}, \Cref{eqn:restriction} and the equality $t_{0} \circ \psi_{g} = \phi_{g} \circ t_{0}$ we deduce that $\psi_{g}^{*}(\sigma'_{t_{0}}) = \sigma'_{t_{0}}$, hence $\sigma'_{t_{0}}$ is a $G_{0}$-invariant adapted differential on $\bar{Y}_{0}'$.
    By \Cref{lm:galoisrestriction} we can find an adapted differential $\sigma_{t}$ on $Y$ such that $(f^{Y}_{0})^{*}(\sigma_{t}) = \sigma'_{t_{0}}$.
    Moreover, the resulting adapted differential is independent of the chosen irreducible component $Y_{0}'$.
    Indeed, let $Y_{1}'$ be another irreducible component and let $\tau_{t}$ be an adapted differential on $Y$ such that $(f^{Y}_{1})^{*}(\tau_{t}) = \sigma'_{t_{1}}$, where $f^{Y}_{1}$ and $t_{1}$ are defined as before but for the irreducible component $Y_{1}'$ instead.
    Since $G$ acts transitively on the irreducible components of the fiber product \cite[Corollary 5.10]{sv96}, we can find some $g \in G$ such that the corresponding automorphism of the fiber product induces an isomorphism $Y_{0}' \to Y_{1}'$ over $Y$, hence also an isomorphism $\psi_{g} \colon \bar{Y}_{0}' \to \bar{Y}_{1}'$ over $Y$.
    From $G$-invariance of $\sigma'$, the equality $\phi_{g} \circ t_{0} = t_{1} \circ \psi_{g}$, the equality $f^{Y}_{0} = f^{Y}_{1} \circ \psi_{g}$ and the usual equations we have
    \[ (f^{Y}_{0})^{*}(\sigma_{t}) = \sigma'_{t_{0}} = \sigma'_{\phi_{g} \circ t_{0}} = \psi_{g}^{*}(\sigma'_{t_{1}}) = \psi_{g}^{*}((f^{Y}_{1})^{*}(\tau_{t})) = (f^{Y}_{0})^{*}(\tau_{t}), \]
    hence $\sigma_{t} = \tau_{t}$ by \Cref{lm:injectiverestriction}.
    
    We show next that the family of adapted differentials $\sigma$ constructed in this manner satisfies \Cref{eqn:limit}, i.e., that $\sigma \in \mathscr{G}(T)$.
    Let $t \colon Y \to T$ be a morphism such that the composition $Y \to X$ is an adapted morphism and let $f \colon Z \to Y$ be a morphism such that the composition $Z \to X$ is an adapted morphism as well.
    Since $f \colon Z \to Y$ is dominant, we may find an irreducible component of $T' \times_{T} Z$ which dominates an irreducible component of $T' \times_{T} Y$, say $Z_{0}' \to Y_{0}'$.
    From the universal property of the normalization we obtain now a morphism $f_{0} \colon \bar{Z}_{0}' \to \bar{Y}_{0}'$ such that $f \circ f_{0}^{Z} = f_{0}^{Y} \circ f_{0}$, where $f_{0}^{Y} \colon \bar{Y}_{0}' \to Y$ and $f_{0}^{Z} \colon \bar{Z}_{0}' \to Z$ are defined as in the last paragraph:
    \begin{center}
      \begin{tikzcd}
        \bar{Z}_{0}' \arrow{r}{f_{0}^{Z}} \arrow[swap]{d}{f_{0}} & Z \arrow{d}{f} \\
        \bar{Y}_{0}' \arrow{r}{f_{0}^{Y}} \arrow[swap]{d}{t_{0}} & Y \arrow{d}{t} \\
        T' \arrow{r}{\pi} & T.
      \end{tikzcd}
    \end{center}
    From \Cref{eqn:limit} for $\sigma'$, the equality $f \circ f_{0}^{Z} = f_{0}^{Y} \circ f_{0}$ and the construction of $\sigma$ we deduce that $(f_{0}^{Z})^{*}(f^{*}(\sigma_{t})) = (f_{0}^{Z})^{*}(\sigma_{t \circ f})$, hence $f^{*}(\sigma_{t}) = \sigma_{t \circ f}$ by \Cref{lm:injectiverestriction}.

    To finish the proof we need to check that $\mathscr{G}(\pi)(\sigma) = \sigma'$, so let $s \colon Y \to T'$ be a morphism such that the composition $Y \to X$ is an adapted morphism.
    Let $f^{Y}_{0} \colon \bar{Y}_{0}' \to Y$ and $t_{0} \colon \bar{Y}_{0}' \to T'$ be constructed as above with respect to the morphism $t := \pi \circ s$.
    For every $g \in G$, denote by $E_{g}$ the equalizer of $\phi_{g} \circ t_{0}$ and $s \circ f_{0}^{Y}$, which is a closed subscheme of $\bar{Y}_{0}'$ because $\pi$ is separated \cite[Definition and Proposition 9.7.(ii)]{gw10}.
    We claim that $\bar{Y}_{0}' = \cup_{g \in G} E_{g}$ as sets.
    Indeed, it suffices to show that every closed point of $\bar{Y}_{0}'$ belongs to $E_{g}$ for some $g \in G$, because all $E_{g}$ are closed subspaces and the set of closed points is dense in $\bar{Y}_{0}'$.
    But we are working over an algebraically closed field, so this follows from \cite[Exercise 9.7]{gw10} and from the equality $\pi \circ t_{0} = \pi \circ s \circ f_{0}^{Y}$, in which $\pi$ is a quotient of $T'$ by $G$.
    Since $\bar{Y}_{0}'$ is irreducible and $G$ is finite, there exists some $g \in G$ such that $\bar{Y}_{0}' = E_{g}$ as topological spaces.
    But $\bar{Y}_{0}'$ is reduced, so we must have $\bar{Y}_{0}' = E_{g}$ and thus $\phi_{g} \circ t_{0} = s \circ f_{0}^{Y}$.
    Let now $\bar{Y}_{1}'$ be a normalized component of the fiber product such that $g$ induces an isomorphism $\psi_{g} \colon \bar{Y}_{0}' \to \bar{Y}_{1}'$ with $t_{1} \circ \psi_{g} = \phi_{g} \circ t_{0}$, notation again as above.
    Then the following diagram remains commutative after adding the dashed arrow:
    \begin{center}
      \begin{tikzcd}
        \bar{Y}_{0}' \arrow[swap]{d}{t_{0}} \arrow{r}{\psi_{g}} \arrow[bend left=45]{rr}{f_{0}^{Y}} & \bar{Y}_{1}' \arrow[swap]{d}{t_{1}} \arrow{r}{f_{1}^{Y}} & Y \arrow{d}{t} \arrow[dashed, swap]{dl}{s} \\
        T' \arrow{r}{\phi_{g}} & T' \arrow{r}{\pi} & T.
      \end{tikzcd}
    \end{center}
    Therefore $t_{1} \circ \psi_{g} = s \circ f^{Y}_{1} \circ \psi_{g}$, and thus also $t_{1} = s \circ f^{Y}_{1}$.
    By construction of $\sigma$ we have $(f^{Y}_{1})^{*}(\sigma_{t}) = \sigma_{t_{1}}'$, so $(f^{Y}_{1})^{*}(\sigma_{t}) = (f^{Y}_{1})^{*}(\sigma'_{s})$.
    It follows from \Cref{lm:injectiverestriction} that $\sigma_{\pi \circ s} = \sigma'_{s}$, hence $\mathscr{G}(\pi)(\sigma)_{s} = \sigma'_{s}$.
  \end{proof}
\end{lm}

\subsection{Conclusion}
\label{subsection:conclusion}
If we denote by $i_{1} \colon \operatorname{Adpt}(X,\Delta) \to \operatorname{Nor}/X$ and $i_{2} \colon \operatorname{Nor}/X \to \operatorname{Sch}/X$ the inclusion functors, \Cref{rem:compositionofextension} implies that
\[ \Omega_{(X,\Delta)}^{p} = {_P}(i_{2} \circ i_{1})(\mathscr{F}) \cong {_P}i_{2}({_P}i_{1}(\mathscr{F})) = {_P}i_{2}(\mathscr{G}). \]
In \Cref{subsection:zariski} and \Cref{subsection:galois} we have seen that $\mathscr{G}$ is a qfh-sheaf on $\operatorname{Nor}/X$ in the sense of Suslin--Voevodsky, i.e., that it satisfies conditions \ref{condition:1} and \ref{condition:2} from the introduction.
Therefore, by \Cref{lm:kan}, its Kan extension $\Omega_{(X,\Delta)}^{p}$ is a qfh-sheaf on $\operatorname{Sch}/X$.

\section{Some computations}
\label{section:computations}

The value $\Omega_{(X,\Delta)}^{p}(T)$ can be described more explicitly for quasi-finite separated $X$-schemes $T \to X$.
In \Cref{prop:globalsections} we compute the restriction of $\Omega_{(X,\Delta)}^{p}$ to the small Zariski site on $X$, which allows us to compute the cohomology of $\Omega_{(X,\Delta)}^{p}$ in \Cref{cor:cohomology}.
The cases $\dim(T) < \dim(X)$ and $\dim(T) = \dim(X)$ are briefly discussed in \Cref{rem:nohomotopyinvariance} and \Cref{rem:qfssections} respectively.

\begin{prop}\label{prop:globalsections}
  In the setting of \Cref{thm:main}, denote by $\Omega_{(X,\Delta)}^{p}|_{X}$ the restriction of $\Omega_{(X,\Delta)}^{p}$ to the small Zariski site on $X$.
  Then we have
  \[ \Omega_{X}^{[p]}(\log{\lfloor \Delta \rfloor}) \cong \Omega_{(X,\Delta)}^{p}|_{X}. \]
  \begin{proof}
    The global sections of $\Omega_{(X,\Delta)}^{p}$ are given by
    \[ \Omega_{(X,\Delta)}^{p}(X) \cong \left\{ (\sigma_{\gamma})_{\gamma} \in \prod_{\substack{\gamma \colon Y \to X \\ \text{adapted}}} \Omega_{(X,\Delta,\gamma)}^{[p]}(Y) \,\middle|\,
    \begin{aligned}
      & f^{*}(\sigma_{\gamma}) = \sigma_{\gamma \circ f} \text{ for all } \\
      & f \colon Y' \to Y \text{ over } X
    \end{aligned}
    \right\}. \]
    Recall that $\Omega_{X}^{[p]}(\log{\lfloor \Delta \rfloor })(X)$ can be described as $i_{*}\Omega_{U}^{p}(\log{ \lfloor \Delta|_{U} \rfloor} )$ for $U \subseteq X$ as in \Cref{defn:adapted}.
    Therefore, given $\sigma \in \Omega_{X}^{[p]}(\log{ \lfloor \Delta \rfloor })(X) = \Omega_{U}^{p}(\log{ \lfloor \Delta|_{U} \rfloor })(U)$ and an adapted morphism $\gamma \colon Y \to X$, we can pull back $\sigma$ along $\gamma|_{\gamma^{-1}(U)} \colon \gamma^{-1}(U) \to U$ and restrict the result to the open subset $V$ in \Cref{defn:adapted} to obtain a section $\gamma^{*}\sigma \in \Omega_{(X,\Delta,\gamma)}^{[p]}(Y)$.
    By torsion freeness and functoriality of pull-back of logarithmic differentials, the resulting family $(\gamma^{*}\sigma)_{\gamma}$ is a section in $\Omega_{(X,\Delta)}^{p}(X)$.
    Since adapted morphisms are dominant, this defines an injective morphism $\Omega_{X}^{[p]}(\log{ \lfloor \Delta \rfloor })(X) \to \Omega_{(X,\Delta)}^{p}(X)$.

    The same recipe gives an injective morphism $\Omega_{X}^{[p]}(\log{\lfloor \Delta \rfloor})(W) \to \Omega_{(X,\Delta)}^{p}(W)$ for any other open subset $W \subseteq X$.
    Suppose now that $\iota \colon W_{1} \to W_{2}$ is an inclusion of open subsets of $X$.
    The restriction morphism of $\Omega_{X}^{[p]}(\log{\lfloor \Delta \rfloor})$ is induced by the pull-back of Kähler differentials $\iota^{*}$ as in \Cref{lm:presheaf}, so we want to show that the following diagram commutes:
    \begin{center}
      \begin{tikzcd}
        \Omega_{X}^{[p]}(\log{\lfloor \Delta \rfloor})(W_{2}) \arrow[swap]{d}{\iota^{*}} \arrow{r} & \Omega_{(X,\Delta)}^{p}(W_{2}) \arrow{d}{\Omega_{(X,\Delta)}^{p}(\iota)} \\
        \Omega_{X}^{[p]}(\log{\lfloor \Delta \rfloor})(W_{1}) \arrow{r} & \Omega_{(X,\Delta)}^{p}(W_{1}).
      \end{tikzcd}
    \end{center}
    Let $\sigma \in \Omega_{X}^{[p]}(\log{\lfloor \Delta \rfloor})(W_{2})$ and let $w_{1} \colon Y \to W_{1}$ be a morphism such that the composition $\gamma_{w_{1}} \colon Y \to X$ is an adapted morphism.
    The upper horizontal isomorphism sends $\sigma$ to the family $(w_{2}^{*}\sigma)_{w_{2}}$, indexed by morphisms $w_{2} \colon Z \to W_{2}$ such that the composition $Z \to X$ is an adapted morphism.
    Its image under $\Omega_{(X,\Delta)}^{p}(\iota)$ is the family whose component at $w_{1}$ is given by $(w_{2}^{*}\sigma)_{w_{2} = \iota \circ w_{1}}$, i.e., $w_{1}^{*}\iota^{*}\sigma$.
    Therefore, the diagram commutes and we have an injective morphism
    \[ \Omega_{X}^{[p]}(\log{\lfloor \Delta \rfloor}) \to \Omega_{(X,\Delta)}^{p}|_{X} \]
    of sheaves on the topological space $X$.
    We show that it is also surjective.
    Suppose we are given a family $(\sigma_{\gamma})_{\gamma} \in \Omega_{(X,\Delta)}^{p}(X)$.
    Let $U_{0} \subseteq X$ be the largest open subset over which $X$ and $\operatorname{Supp}(\Delta)$ are smooth, which is a big open subset contained in the open subset $U \subseteq X$ from \Cref{defn:adapted}.
    Surjectivity is local on $X$, so after possibly shrinking $X$, we may assume that we are in the setting of \cite[\S 3.5]{ev92}.
    Therefore, we may assume that there exists an adapted morphism $\gamma_{0} \colon Y_{0} \to X$ factoring through $U_{0} \subseteq X$ such that there exists a finite group $G$ acting on $Y_{0}$ so that $Y_{0} \to U_{0}$ is a quotient of $Y_{0}$ by $G$.
    \Cref{eqn:limit} ensures that the logarithmic differential $\sigma_{\gamma_{0}}$ is $G$-invariant, hence Hurwitz's formula \cite[Lemma 3.16]{ev92} implies that there exists a logarithmic differential $\sigma \in \Omega_{U_{0}}^{p}(\log{ \lfloor \Delta|_{U_{0}} \rfloor })(U_{0})$ such that $\gamma_{0}^{*}\sigma = \sigma_{\gamma_{0}}$.
    This logarithmic differential defines a reflexive logarithmic differential in $\Omega_{X}^{[p]}(\log{ \lfloor \Delta \rfloor })(X)$ which is the desired preimage of the given section of $\Omega_{(X,\Delta)}^{p}$.
    The same argument works for sections defined over smaller open subsets $W \subseteq X$, so the morphism above is surjective.
  \end{proof}
\end{prop}

\begin{cor}
  \label{cor:cohomology}
  In the setting of \Cref{thm:main} we have
  \[ H^{i}_{\mathrm{qfh}}\left(X,\Omega_{(X,\Delta)}^{p}\right) \cong H^{i}\left(X,\Omega_{X}^{[p]}(\log{\lfloor \Delta \rfloor})\right) \]
  for all $i \in \mathbb{N}$, where the left-hand side denotes cohomology with respect to the qfh-topology on $\operatorname{Sch}/X$ and the right-hand side denotes usual sheaf cohomology on $X$.
  \begin{proof}
    Since normal schemes are geometrically unibranch \cite[\href{https://stacks.math.columbia.edu/tag/0BQ3}{0BQ3}]{stacks-project} and qfh-sheaves are étale sheaves, we can apply \cite[Proposition 10.5.12]{cd19} and \cite[Theorem 10.5.10]{cd19} to compute the left-hand side with respect to the étale topology on $\operatorname{Sch}/X$.
    Étale cohomology groups can be computed in the small étale site instead \cite[\href{https://stacks.math.columbia.edu/tag/03YX}{03YX}]{stacks-project}.
    By \Cref{prop:globalsections}, $\Omega_{(X,\Delta)}^{p}|_{X}$ is a quasi-coherent sheaf on $X$, so we can in turn compute its étale cohomology groups as usual Zariski sheaf cohomology on $X$ \cite[\href{https://stacks.math.columbia.edu/tag/03DW}{03DW}]{stacks-project}.
  \end{proof}
\end{cor}

\begin{rem}\label{rem:nohomotopyinvariance}
  In the setting of \Cref{thm:main}, let $T \to X$ be an $X$-scheme such that $\dim(T) < \dim(X)$.
  Then $\Omega_{(X,\Delta)}^{p}(T) = 0$, because we are taking the limit over an empty diagram in \Cref{eqn:adjoint}.
  But \Cref{prop:globalsections} shows that $\Omega_{(X,\Delta)}^{p}(X) \neq 0$ for some $(X,\Delta)$, so $\Omega_{(X,\Delta)}^{p}$ is not homotopy invariant in general.
\end{rem}

\begin{rem}\label{rem:qfssections}
  If $T \to X$ is quasi-finite with $\dim(T) = \dim(X)$, then the normalization $\bar{T} \to T$ induces a canonical isomorphism $\Omega_{(X,\Delta)}^{p}(T) \cong \Omega_{(X,\Delta)}^{p}(\bar{T})$ as a consequence of its universal property \cite[\href{https://stacks.math.columbia.edu/tag/035Q}{035Q}]{stacks-project}.
  Combining this with the Zariski sheaf property one can express $\Omega_{(X,\Delta)}^{p}(T)$ as the product of the sections over its irreducible components, and with arguments similar to those in \Cref{prop:globalsections} this allows one to compute the sections of $\Omega_{(X,\Delta)}^{p}$ over any quasi-finite separated $X$-scheme $T \to X$.
\end{rem}

\printbibliography
\vfill

\end{document}